\documentclass[a4paper]{amsart}
\usepackage{amssymb,enumerate}
\usepackage{epsfig}

\title{Interval orders and reverse mathematics}
\author{Alberto Marcone}
    \address{Dipartimento di Matematica e Informatica,
    Universit\`{a} di Udine,
    viale delle Scienze 206,
    33100 Udine,
    Italy}
\email{marcone@dimi.uniud.it}


\date{December 16, 2006}

\theoremstyle{plain}
\newtheorem{theorem}{Theorem}[section]
\newtheorem{characterization}{Characterization Theorem}
\newtheorem{lemma}[theorem]{Lemma}
\newtheorem{cor}[theorem]{Corollary}
\newtheorem{claim}{Claim}

\theoremstyle{definition}
\newtheorem{definition}[theorem]{Definition}
\newtheorem{remark}[theorem]{Remark}

\newcommand{\set}[2]{\left\{\,{#1}\mid{#2}\,\right\}}
\newcommand{\N}{{\mathbb N}}
\newcommand{\Z}{{\mathbb Z}}
\newcommand{\lh}{\operatorname{lh}}
\newcommand{\PI}[2]{\ensuremath{\boldsymbol\Pi^{#1}_{#2}}}
\newcommand{\SI}[2]{\ensuremath{\boldsymbol\Sigma^{#1}_{#2}}}
\newcommand{\DE}[2]{\ensuremath{\boldsymbol\Delta^{#1}_{#2}}}
\newcommand{\system}[1]{\mbox{\fontfamily{cmss}\fontshape{n}\fontseries{m}%
    \selectfont#1}}

\newcommand{\RCA}{\system{RCA}\ensuremath{_0}}
\newcommand{\WKL}{\system{WKL}\ensuremath{_0}}
\newcommand{\ACA}{\system{ACA}\ensuremath{_0}}

\newcommand{\PP}{\ensuremath{\mathbf{P}}}
\newcommand{\LL}{\ensuremath{\mathbf{L}}}
\newcommand{\tpt}{\ensuremath{\mathbf{2} \oplus \mathbf{2}}}
\newcommand{\tpo}{\ensuremath{\mathbf{3} \oplus \mathbf{1}}}
\newcommand{\dw}[1]{\ensuremath{{{#1} \!\downarrow}}}
\newcommand{\uw}[1]{\ensuremath{{{#1} \!\uparrow}}}

\begin{document}
\pagestyle{plain}
\begin{abstract}
We study the reverse mathematics of interval orders. We establish
the logical strength of the implications between various definitions
of the notion of interval order. We also consider the strength of
different versions of the characterization theorem for interval
orders: a partial order is an interval order if and only if it does
not contain \tpt. We also study proper interval orders and their
characterization theorem: a partial order is a proper interval order
if and only if it contains neither \tpt\ nor \tpo.
\end{abstract}
\maketitle

Interval orders are a particular kind of partial orders which occur
quite naturally in many different areas and have been widely
studied. A partial order $\PP = (P, {\leq_P})$ is an \emph{interval
order} if the elements of $P$ can be mapped to nonempty intervals of
a linear order \LL\ so that $p <_P q$ holds iff every element of the
interval associated to $p$ precedes every element of the interval
associated to $q$. The linear order \LL\ and the map from $P$ to
intervals are called an \emph{interval representation of \PP}. The
basic reference on interval orders is Fishburn's monograph
\cite{Fishburn}.

The name \lq\lq interval order\rq\rq\ was introduced by Fishburn
(\cite{Fishburn-char}), although the notion was already studied much
earlier by Norbert Wiener (\cite{Wiener}), who used the terminology
\lq\lq relation of complete sequence\rq\rq. Interval orders model
many phenomena occurring in the applied sciences:
\cite[\S2.1]{Fishburn} include examples such as chronological dating
in archaeology and paleontology, scheduling of manufacturing
processes, and psychophysical perception of sounds. Notice that if
\PP\ is a countable interval order then we can assume that \LL\ is
the rational or (as usual in applications) the real line (a
\emph{real representation}, in the terminology of \cite{Fishburn}).

Most recent research on interval orders (see e.g.\ the survey
\cite{Trotter} and chapter 8 of \cite{schroder}) focuses on finite
partial orders, while in this paper we consider mostly infinite ones
(although a careful analysis of the finite case is instrumental in
obtaining results in the infinite case). A recent result about
infinite interval orders shows that every interval order which is a
well quasi-order is a better quasi-order (\cite{Pouzet-Sauer}).

The basic characterization for interval orders is given by the
following theorem proved independently by Fishburn
(\cite{Fishburn-char}) and Mirkin (\cite{Mirkin2}):

\begin{characterization}\label{char}
A partial order is an interval order if and only if it does not
contain \tpt.
\end{characterization}

Here \lq\lq\PP\ does not contain \tpt\rq\rq\ means that for no $P'
\subseteq P$ the restriction of $\leq_P$ to $P'$ is the partial
order with Hasse diagram\ \
\raisebox{-2pt}{%
\begin{picture}(10,10) \put(0,0){\circle*{3}}
\put(0,10){\circle*{3}} \put(10,0){\circle*{3}}
\put(10,10){\circle*{3}} \put(0,0){\line(0,1){10}}
\put(10,0){\line(0,1){10}}
\end{picture}
}. It is easy to see that \PP\ does not contain \tpt\ if and only if
\[
\forall p_0, q_0, p_1, q_1 \in P (p_0 \leq_P q_0 \land p_1 \leq_P
q_1 \implies p_0 \leq_P q_1 \lor p_1 \leq_P q_0).
\]

Two natural ways of strengthening the notion of interval order lead
to the definitions of unit interval order and proper interval order.

An interval order with a real representation such that all intervals
have the same positive length (which can be assumed to be $1$) is
called a \emph{unit interval order}.

If an interval order \PP\ has an interval representation such that
an interval associated to an element of $P$ is never a proper subset
of another such interval, then we say that \PP\ is a \emph{proper
interval order}. An interval representation with the above property
is called a \emph{proper interval representation}.

It is immediate that every unit interval order is a proper interval
order. If the partial order is finite then the reverse implications
is also true (\cite{Roberts}, see \cite{Bogart-West} for a short
proof). On the other hand, there exist infinite proper interval
orders which are not unit interval orders: a simple example is
provided by the ordinal $\omega+1$. Notice however that the fact
that $\omega+1$ is not a unit interval order has more to do with the
real line (which in this context appears to be \lq\lq too short'')
than with structural properties of the partial order. Therefore when
dealing with infinite partial orders the notion of proper interval
order appears to be more natural, as witnessed also by the following
characterization theorem.

\begin{characterization}\label{char-proper}
A partial order is a proper interval order if and only if it
contains neither \tpt\ nor \tpo.
\end{characterization}

\lq\lq\PP\ does not contain \tpo\rq\rq\ means that for no $P'
\subseteq P$ the restriction of $\leq_P$ to $P'$ is the partial
order with Hasse diagram\ \
\raisebox{-2pt}{%
\begin{picture}(10,10) \put(0,0){\circle*{3}}
\put(0,5){\circle*{3}} \put(0,10){\circle*{3}}
\put(10,5){\circle*{3}} \put(0,0){\line(0,1){10}}
\end{picture}
}. It is easy to see that \PP\ does not contain \tpo\ if and only if
\[
\forall p_0, p_1, p_2, q \in P (p_0 <_P p_1 <_P p_2 \implies p_0
\leq_P q \lor q \leq_P p_2).
\]

Characterization Theorem \ref{char-proper} is usually known as the
Scott-Suppes Theorem. Scott and Suppes (\cite{Scott-Suppes}) proved
the theorem in the finite case for unit interval orders (see
\cite{Balof-Bogart} for a simple proof in this setting). Fishburn's
monograph includes a proof of this theorem with no restrictions on
cardinality (\cite[Theorem 2.7]{Fishburn}).

\bigskip

In this paper we study interval orders and proper interval orders
from the viewpoint of \emph{reverse mathematics}. The basic
reference for reverse mathematics is Simpson's book \cite{sosoa},
which contains all background material needed for this paper (and
much more). A sample of recent research in the area is contained in
\cite{RM2001}.

In reverse mathematics, one formalizes theorems of ordinary
mathematics and attempts to discover the set theoretic axioms
required to prove these theorems. This project is usually carried
out in the context of subsystems of second order arithmetic, taking
\RCA\ as the base system. \RCA\ is the subsystem obtained from full
second order arithmetic by restricting the comprehension scheme to
\DE01 formulas and adding a formula induction scheme for \SI01
formulas. In this paper, we will be concerned only with \RCA\ and
its fairly weak extension known as \WKL\ (\WKL\ is strictly weaker
than the subsystem \ACA\ obtained by extending the comprehension
scheme in \RCA\ to all arithmetic formulas). \WKL\ is obtained by
adjoining Weak K\"{o}nig's Lemma (i.e.\ K\"{o}nig's Lemma for trees
of sequences of $0$'s and $1$'s) to \RCA.

Many results about partial and linear orders have been studied from
the viewpoint of reverse mathematics: recent papers include
\cite{Dushnik-Miller,DHLS,CMS,CR,Friedcomp,hirst-survey,wqobqo,Montalban}.
Moreover, \cite[\S3]{Schmerl} includes a couple of results about
interval graphs, which are strictly connected to interval orders.

\section{Overview of results and plan of the paper}
The first step in the study of a new topic in the context of reverse
mathematics is finding appropriate formalizations of the relevant
notions. Often, this requires making choices between classically
equivalent definitions for the mathematical concepts appearing in
the definitions. In this paper, we consider a number of equivalent
definitions for the notions of interval order and of proper interval
order, and we examine how difficult it is to prove the equivalences
of these definitions.

There is no particular difficulty in coding a countable partial
order in the weak base theory \RCA. The only point to note is that
we consider only countable partial orders.

However the notion of interval order hinges on the notion of
interval of a linear order, and the latter can be interpreted in
different ways, leading to notions that are not necessarily
equivalent in the weak base theory \RCA. We can define an interval
of the linear order $\LL = (L, {\leq_L})$ to be a set $I \subseteq
L$ which satisfies $\forall x,y \in I\, \forall z \in L (x \leq_L z
\leq_L y \implies z \in I)$. Another possibility is to restrict our
attention to closed intervals (this is often done in the literature
about interval orders, e.g.\ in \cite{Trotter} this is done from the
outset) and code them by pairs $(a,b)$ of elements of $L$ such that
$a \leq_L b$ (obviously in this case $x \in L$ belongs to the
interval if and only if $a \leq_L x \leq_L b$). If we apply the
latter concept of interval we speak of a \emph{closed interval
representation} of the partial order. In defining interval orders
there is a further subtlety, that turns out to be important in our
study of the proof theoretic strength of various statements: i.e.\
we may require the map of the interval representation to be
injective. Combining the two possible choices in each of the two
cases we obtain four notions of interval order: interval order, 1-1
interval order, closed interval order, and 1-1 closed interval
order. Another notion is obtained by further strengthening the
definition of 1-1 closed interval order: a closed interval
representation is a \emph{distinguishing representation} if all
endpoints of the closed intervals are distinct (see e.g.\
\cite{Trotter}). This leads to the notion of distinguishing interval
order. In Section \ref{sect:def} we will give the precise
definitions of these notions in \RCA.

The five notions introduced above are all equivalent, and we
establish the axioms needed to show the equivalences among them and
with the characterization provided by Characterization Theorem
\ref{char}. (Notice that the proofs of the latter theorem in
\cite{Fishburn} and \cite{Trotter} can be easily carried out in
\ACA: see Remark \ref{ACA} below.)

\begin{figure}
\includegraphics{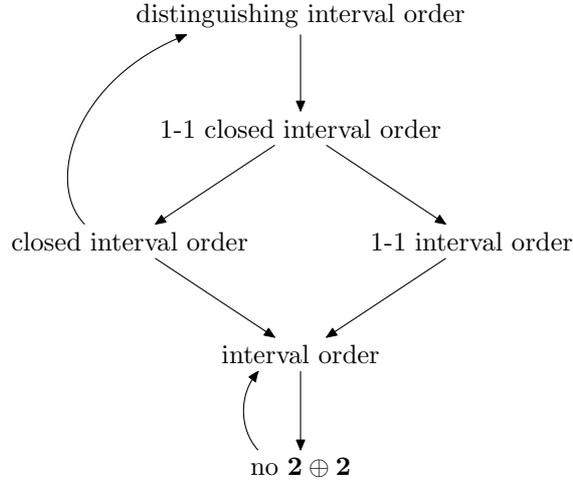}
\caption{\label{diagram}Implications about interval orders provable
in \RCA.}
\end{figure}

We show that \RCA\ proves exactly the implications appearing in
Figure \ref{diagram} (where an arrow with origin in the node labeled
$A$ pointing towards the node labeled $B$ represents the statement
\lq\lq every partial order which satisfies $A$ satisfies $B$\rq\rq),
or that can be obtained by composing arrows appearing in that
diagram.
In particular we obtain the following result about Characterization
Theorem \ref{char}:

\begin{theorem}\label{BCT-RCA}
\RCA\ proves that a partial order is an interval order if and only
if it does not contain \tpt.
\end{theorem}

The arrows pointing downwards (possibly diagonally) in Figure
\ref{diagram} either follow from the definitions or are
straightforward to prove (these implications are collected in
Theorem \ref{obvious}), while the two arrows pointing upwards will
be proved in \S\ref{sect:RCA}.

Figure \ref{diagram} implies that in \RCA\ there are at most three
distinct notions of interval order. In order of decreasing strength
these are: closed interval order, 1-1 interval order, and interval
order. In Section \ref{sect:WKL} we show that each of the missing
implications is equivalent to \WKL. For the stronger notions of
interval order we obtain the following reverse mathematics results
about Characterization Theorem \ref{char}:

\begin{theorem}
In \RCA\ the following are equivalent:
\begin{enumerate}[(i)]
 \item
\WKL;
 \item
a partial order is a 1-1 interval order if and only if it does not
contain \tpt;
 \item
a partial order is a closed interval order if and only if it does
not contain \tpt;
 \item
a partial order is a 1-1 closed interval order if and only if it
does not contain \tpt;
 \item
a partial order is a distinguishing interval order if and only if it
does not contain \tpt.
\end{enumerate}\end{theorem}

In particular this implies that \RCA\ does not prove that the
equivalence between the three notions of interval order mentioned
above.

Section \ref{sect:finite} is devoted to a detailed analysis of the
equivalences for finite partial orders; this analysis will be used
in the proofs of the following sections.

\medskip

\begin{figure}
\includegraphics{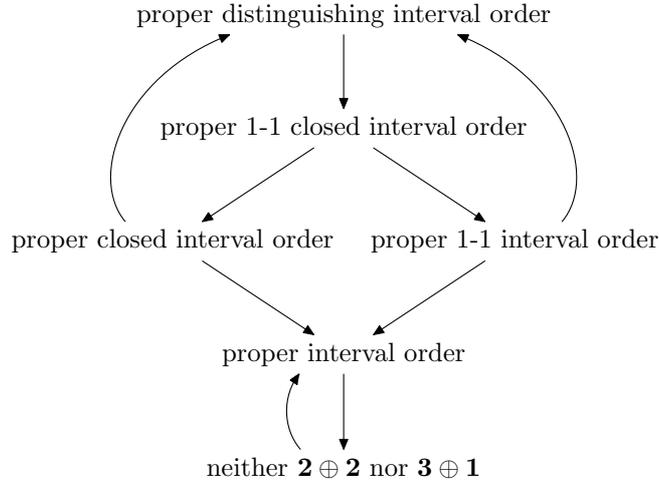}
\caption{\label{diagram-proper}Implications about proper interval
orders provable in \RCA.}
\end{figure}

When defining proper interval orders the same choices about
intervals and injectivity are possible: we thus also have five
different notions of proper interval order, plus the
characterization provided by Characterization Theorem
\ref{char-proper}. We show that \RCA\ proves exactly the
implications appearing in Figure \ref{diagram-proper}, or that can
be obtained by composing arrows appearing in that diagram.
In particular we obtain the following result about Characterization
Theorem \ref{char-proper}:

\begin{theorem}\label{BCT-RCA-proper}
\RCA\ proves that a partial order is a proper interval order if and
only if it contains neither \tpt\ nor \tpo.
\end{theorem}

Figures \ref{diagram} and \ref{diagram-proper} are similar, except
that the latter includes one arrow whose analogous is missing from
the former. Indeed within \RCA, a 1-1 interval order is necessarily
a distinguishing interval order if we have a proper representation,
but not in general.

Figure \ref{diagram-proper} implies that in \RCA\ there are at most
two distinct notions of proper interval order, i.e.\ proper closed
interval order and proper interval order. We show that the missing
implication is equivalent to \WKL, even if we restrict ourselves to
closed interval orders. For the stronger notions of interval order
we obtain the following reverse mathematics results about
Characterization Theorem \ref{char-proper}:

\begin{theorem}
In \RCA\ the following are equivalent:
\begin{enumerate}[(i)]
 \item
\WKL;
 \item
a partial order is a proper 1-1 interval order if and only if it
contains neither \tpt\ nor \tpo;
 \item
a partial order is a proper closed interval order if and only if it
contains neither \tpt\ nor \tpo;
 \item
a partial order is a proper 1-1 closed interval order if and only if
it contains neither \tpt\ nor \tpo;
 \item
a partial order is a proper distinguishing interval order if and
only if it contains neither \tpt\ nor \tpo.
\end{enumerate}\end{theorem}

In Section \ref{sect:proper} the definitions and the arguments of
Sections \ref{sect:def} through \ref{sect:WKL} are adapted to the
case of proper interval orders, and all results about proper
interval orders are proved. Some of the proofs are straightforward
translations of the corresponding proofs for interval orders, while
others exploit the properties of proper interval orders.

\bigskip

Our results are stated in terms of subsystems of second order
arithmetic, but have corollaries that can be viewed as examples of
computable mathematics in the style of \cite{HandRecMath}. Samples
of these corollaries are the following, where we use standard
terminology from computability theory:

\begin{cor}
For every computable partial order \PP\ not containing \tpt\ there
exist a computable linear order \LL\ and a computable function from
$P$ to intervals of \LL\ witnessing that \PP\ is an interval order.
\end{cor}

\begin{cor}
There exists a computable partial order \PP\ not containing \tpt\
such that for every computable linear order \LL\ there is no
computable function from $P$ to closed intervals of \LL\ witnessing
that \PP\ is a closed interval order.
\end{cor}

\begin{cor}
For every computable partial order \PP\ not containing \tpt\ there
exist a low (resp.\ almost recursive) linear order \LL\ and a low
(resp.\ almost recursive) function from $P$ to closed intervals of
\LL\ witnessing that \PP\ is a distinguishing interval order.
\end{cor}

(The last Corollary follows from our results by the properties of
$\omega$-models of \WKL\ which appear in \cite[\S VIII.2]{sosoa}.)
\bigskip

We assume some familiarity of the reader with subsystems of second
order arithmetic, but the paper is self-contained as far as interval
order theory is concerned.

From now on, when a definition or the statement of a result starts
with the name of a subsystem of second order arithmetic in
parenthesis, it means that the definition is given, or the statement
provable, in that subsystem.

\section{Definitions and elementary facts}\label{sect:def}
\begin{definition}
(\RCA) A \emph{partial order} \PP\ is a pair $(P, {\leq_P})$ where
$P$ is a set and ${\leq_P} \subseteq P \times P$ is reflexive,
transitive and anti-symmetric. The partial order \PP\ is a
\emph{linear order} if we have also $\forall p,q \in P (p \leq_P q
\lor q \leq_P p)$.
\end{definition}

\begin{remark}
If \PP\ is a partial order then $P \subseteq \N$ and hence on $P$ we
have also the restriction of the usual order on the natural numbers.
When there is danger of confusion we denote the latter by $\leq_\N$.
\end{remark}

\begin{definition}
(\RCA) If \PP\ is a partial order we define the relations ${<_P}$
and ${\perp_P}$ as follows:
\begin{gather*}
p <_P q \iff p \leq_P q \land p \neq q,\\
p \perp_P q \iff p \nleq_P q \land q \nleq_P p.
\end{gather*}
\end{definition}

Sometimes it is convenient to use quasi-orders, which are defined by
dropping the requirement of anti-symmetry from the definition of
partial order. In particular we will be interested in linear
quasi-orders.

\begin{definition}
(\RCA) $\PP = (P, {\leq_P})$ is a \emph{quasi-order} if ${\leq_P}
\subseteq P \times P$ is reflexive and transitive. If we have also
$\forall p,q \in P (p \leq_P q \lor q \leq_P p)$ we say that \PP\
is a \emph{linear quasi-order}.
\end{definition}

\begin{definition}
(\RCA) If \PP\ is a quasi-order we define ${<_P}$ by
\[
p <_P q \iff p \leq_P q \land p \nleq_P q,
\]
while no changes are needed in the definition of $\perp_P$.
Furthermore we define $\equiv_P$ by
\[
p \equiv_P q \iff p \leq_P q \land p \leq_P q.
\]
\end{definition}

It is immediate to check in \RCA\ that if \PP\ is a quasi-order then
$\equiv_P$ is an equivalence relation.

In our setting using (linear) quasi-orders in place of partial
(resp.\ linear) orders is just a matter of convenience, as the
following easy lemma shows.

\begin{lemma}\label{qo}
(\RCA) Let \PP\ be a quasi-order. Then there exist $P' \subseteq P$
and $f:P \to P'$ such that $\PP' = (P', {\leq_P})$ is a partial
order and $f$ is a surjective order-preserving function satisfying
$f(p)=p$ for every $p \in P'$.

Furthermore, if \PP\ is a linear quasi-order then $\PP'$ is a linear
order.
\end{lemma}
\begin{proof}
Since $P \subseteq \N$ we can let
\begin{gather*}
P' = \set{p \in P}{\forall q <_\N p\; q \not\equiv_P p};\\
f(p) = \text{the $<_\N$-least $q$ such that } q \equiv_P p.\qedhere
\end{gather*}
\end{proof}

We can now introduce the different notions of interval order.

\begin{definition}\label{def:io}
(\RCA) A partial order \PP\ is an \emph{interval order} if there
exist a linear order \LL\ and a set $F \subseteq P \times L$ such
that, abbreviating $\set{x \in L}{(p,x) \in F}$ by $F(p)$ for every
$p \in P$, we have:\begin{enumerate}[\quad({i}1)]
 \item
$F(p) \neq \emptyset$ and $\forall x,y \in F(p)\, \forall z \in L (x
<_L z <_L y \implies z \in F(p))$ for all $p \in P$;
 \item
$p <_P q \iff \forall x \in F(p)\, \forall y \in F(q)\, x <_L y$ for
all $p, q \in P$.
\end{enumerate}
\PP\ is a \emph{1-1 interval order} if we have
also\begin{enumerate}[\quad({i}1)]\setcounter{enumi}2
 \item
$F(p) \neq F(q)$ whenever $p \neq q$.
\end{enumerate}

\PP\ is a \emph{closed interval order} if there exist a linear order
\LL\ and two functions $f_0,f_1: P \to L$ such
that:\begin{enumerate}[\quad(c1)]
 \item
$f_0(p) \leq_L f_1 (p)$ for all $p \in P$;
 \item
$p <_P q \iff f_1(p) <_L f_0(q)$ for all $p, q \in P$.
\end{enumerate}
\PP\ is a \emph{1-1 closed interval order} if we have
also\begin{enumerate}[\quad(c1)]\setcounter{enumi}2
 \item
$f_0(p) \neq f_0(q)$ or $f_1(p) \neq f_1(q)$ whenever $p \neq q$.
\end{enumerate}
\PP\ is a \emph{distinguishing interval order} if beside (c1--2) we
have also\begin{enumerate}[\quad(c1)]\setcounter{enumi}3
 \item
$f_i(p) \neq f_j(q)$ whenever $p \neq q$ or $i \neq j$.
\end{enumerate}\end{definition}

It is immediate that if we set $F(p) = \set{x \in L}{f_0(p) \leq_L x
\leq_L f_1(p)}$, conditions (c1--3) are the translations of
conditions (i1--3).

\begin{remark}
Lemma \ref{qo} implies that in the preceding definitions we can use
linear quasi-orders in place of linear orders. Whenever it is
convenient for the clarity of the exposition, we will use this fact
without mentioning it explicitly.
\end{remark}

\begin{definition}
(\RCA) A partial order \PP\ \emph{does not contain \tpt} if
\[
\forall p_0, q_0, p_1, q_1 \in P (p_0 <_P q_0 \land p_1 <_P q_1
\implies p_0 \leq_P q_1 \lor p_1 \leq_P q_0).
\]
\end{definition}

\begin{definition}
(\RCA) If \PP\ is a partial order and $p \in P$ the \emph{strict
downward and upward closures of $p$ in $P$} are the sets
\[
\uw p^\PP = \set{q \in P}{p <_P q} \quad \text{and} \quad \dw p^\PP
= \set{q \in P}{q <_P p}.
\]
When \PP\ is clear from the context we write $\uw p$ and $\dw p$.
\end{definition}

The next lemma is a basic observation about partial orders not
containing \tpt.

\begin{lemma}\label{comparability}
(\RCA) If \PP\ does not contain \tpt\ then for every $p,q \in P$ we
have either $\uw p \subseteq \uw q$ or $\uw q \subseteq \uw p$, and
similarly either $\dw p \subseteq \dw q$ or $\dw q \subseteq \dw p$.
\end{lemma}
\begin{proof}
If $\uw p \nsubseteq \uw q$ and $\uw q \nsubseteq \uw p$ let $p_1
\in \uw p \setminus \uw q$ and $q_1 \in \uw q \setminus \uw p$. Then
$p, p_1, q, q_1$ show that \PP\ contains \tpt.

The argument for the strict downward closures is similar.
\end{proof}

The following lemma is useful to show that an interval order is
actually a 1-1 interval order.

\begin{lemma}\label{1-1condition}
Suppose \PP\ is an interval order such that
\[
\forall p,q \in P (p \neq q \implies \uw p \neq \uw q \lor \dw p
\neq \dw q).
\]
Then \PP\ is a 1-1 interval order.
\end{lemma}
\begin{proof}
Let \LL\ and $F$ satisfy conditions (i1--2). We claim that $F$
satisfies also (i3). Fix $p,q \in P$ with $p \neq q$. We have either
$\uw p \neq \uw q$ or $\dw p \neq \dw q$. Without loss of
generality, we may assume the former inequality holds and there
exists $r \in \uw p \setminus \uw q$. Then $q \nless_P r$ and for
some $x \in F(r)$ and $y \in F(q)$ we have $x \leq_L y$. On the
other hand $p <_P r$ so that $z <_L x$ for all $z \in F(p)$. Hence
$y \notin F(p)$ and $F(p) \neq F(q)$.
\end{proof}

We now prove the \lq\lq easy\rq\rq\ arrows appearing in Figure
\ref{diagram}.

\begin{theorem}\label{obvious}
(\RCA)\begin{enumerate}[(i)]
 \item
Every distinguishing interval order is a 1-1 closed interval order.
 \item
Every 1-1 (closed) interval order is a (closed) interval order.
 \item
Every (1-1) closed interval order is a (1-1) interval order.
 \item
Every interval order does not contain \tpt.
\end{enumerate}\end{theorem}
\begin{proof}
The statements in (i) and (ii) follow immediately from the
definitions (since condition (c4) implies condition (c3)).

For the statements in (iii), given \LL, $f_0$ and $f_1$ as in the
definition of closed interval order let
\[
F = \set{(p,x) \in P \times L}{f_0(p) \leq_L x \leq_L f_1(p)}.
\]

To prove (iv), let $L$ and $F$ witness that \PP\ is an interval
order. Suppose towards a contradiction that $p_0, q_0, p_1, q_1 \in
P$ are such that $p_0 <_P q_0$, $p_1 <_P q_1$, $p_0 \nleq_P q_1$ and
$p_1 \nleq_P q_0$. The third condition implies the existence of $x,y
\in L$ such that $x \in F(p_0)$, $y \in F(q_1)$, and $y \leq_L x$.
Similarly by the fourth condition there exist $x',y'$ such that $x'
\in F(p_1)$, $y' \in F(q_0)$, and $y' \leq_L x'$. The first two
conditions imply respectively $x <_L y'$ and $x' <_L y$: using
transitivity we have $x <_L x$, which is impossible.
\end{proof}

\section{Finite interval orders}\label{sect:finite}
We start by introducing one of the basic tools in the analysis of
partial orders not containing \tpt. Within \RCA\ we can define it
only for finite partial orders.

\begin{definition}\label{def:conjoint}
(\RCA) Given a finite partial order \PP, let $P^+ = \set{p^+}{p \in
P}$, $P^- = \set{p^-}{p \in P}$, and $P^* = P^+ \cup P^-$. Define a
binary relation ${\leq_\PP^*}$ on $P^*$ as follows:
\begin{align*}
p^+ \leq_\PP^* q^+ & \iff \uw p^\PP \supseteq \uw q^\PP;\\
p^- \leq_\PP^* q^- & \iff \dw p^\PP \subseteq \dw q^\PP;\\
p^+ \leq_\PP^* q^- & \iff p <_P q;\\
p^- \leq_\PP^* q^+ & \iff q \nless_P p.
\end{align*}
$\PP^* = (P^*, {\leq_\PP^*})$ is the \emph{conjoint linear
quasi-order associated to \PP}. When \PP\ is clear from the context
we write $\leq^*$ in place of $\leq_\PP^*$.
\end{definition}

The following lemma justifies the use of the words \lq\lq linear
quasi-order\rq\rq\ in Definition \ref{def:conjoint}.

\begin{lemma}\label{finite cio}
(\RCA) If \PP\ is a finite partial order which does not contain
\tpt\ then ${\leq^*}$ is a linear quasi-order.

Moreover $\PP^*$ and the functions $p \mapsto p^-$, $p \mapsto p^+$
show that \PP\ is a closed interval order.
\end{lemma}
\begin{proof}
Reflexivity of ${\leq^*}$ follows immediately from the definition.
Using Lemma \ref{comparability} it is also immediate that for every
$x,y \in P^*$ we have $x \leq^* y$ or $y \leq^* x$.

It remains to show that ${\leq^*}$ is transitive and to this end we
need to consider eight cases. We tackle three of them, the others
being trivial or similar to one of these:\begin{itemize}
 \item
if $p^+ \leq^* q^+ \leq^* r^-$ then $\uw p \supseteq \uw q$ and $q
<_P r$, i.e.\ $r \in \uw q$; therefore $r \in \uw p$ which means $p
<_P r$ and hence $p^+ \leq^* r^-$;
 \item
if $p^+ \leq^* q^- \leq^* r^+$ then $p <_P q$ and $r \nless_P q$;
hence $q \in \uw p \setminus \uw r$ and, by Lemma
\ref{comparability}, $\uw p \supset \uw r$ holds, so that $p^+
\leq^* r^+$;
 \item
if $p^+ \leq^* q^- \leq^* r^-$ then $p <_P q$ and $\dw q \subseteq
\dw r$, which imply $p <_P r$ and hence $p^+ \leq^* r^-$.
\end{itemize}

Since for every $p$ we have $p^- \leq^* p^+$ (in fact $p^- <^* p^+$)
condition (c1) of Definition \ref{def:io} is satisfied. Condition
(c2) follows immediately from the definition.
\end{proof}

\begin{remark}\label{all1type}
Notice that for all $p,q \in P$ we have $p^+ \not\equiv^* q^-$. In
other words, each $\equiv^*$-equivalence class is contained in
either $P^+$ or $P^-$.
\end{remark}

Lemma \ref{finite cio} does not prove that \PP\ is a distinguishing
interval order, or even a 1-1 closed interval order: if $p,q \in P$
are distinct and such that $\dw p = \dw q$ and $\uw p = \uw q$ we
have $p^- \equiv^* q^-$ and $p^+ \equiv^* q^+$. To obtain the
stronger conclusions we can proceed as follows.

\begin{definition}\label{compatible}
(\RCA) Given a finite partial order \PP\ which does not contain
\tpt, let $\PP^*$ be the conjoint linear quasi-order associated to
\PP. A linear order $(P^*, \leq_L)$ is \emph{compatible with
$\PP^*$} if
\[
\forall x,y \in P^* (x <^* y \implies x <_L y).
\]
\end{definition}

\begin{remark}
Each $(P^*, \leq_L)$ compatible with $\PP^*$ is defined by giving a
linear order on each ${\equiv^*}$-equivalence class, and keeping the
order between ${\equiv^*}$-inequivalent elements unchanged.
\end{remark}

\begin{lemma}\label{existcompatible}
(\RCA) If \PP\ is a finite partial order which does not contain
\tpt\ then there exists a linear order compatible with $\PP^*$.
\end{lemma}
\begin{proof}
For example let
\[
x \leq_L y \iff x <^* y \lor (x \equiv^* y \land x \leq_\N y).
\]
$\leq_L$ is a linear order compatible with $\PP^*$.
\end{proof}

\begin{lemma}\label{finite BCT}
(\RCA) Any finite partial order which does not contain \tpt\ is a
distinguishing interval order.
\end{lemma}
\begin{proof}
Let \PP\ be a finite partial order which does not contain \tpt, and,
by Lemma \ref{existcompatible}, $\leq_L$ a linear order compatible
with $\PP^*$. By Lemma \ref{finite cio} and Remark \ref{all1type}
$(P^*, {\leq_L})$ and the functions $p \mapsto p^-$, $p \mapsto p^+$
witness that \PP\ is a distinguishing interval order.
\end{proof}

Combining Lemma \ref{finite BCT} with Theorem \ref{obvious} we
obtain that \RCA\ proves the equivalence of the six
characterizations of interval orders restricted in the case of
finite partial orders.

\begin{remark}\label{ACA}
The reader should notice that we carried out the discussion in this
section only for finite partial orders, but the constructions and
arguments apply also for infinite ones. However in the infinite case
\RCA\ does not suffice to define $\leq^*$ and we need to use \ACA.
Indeed, arithmetical comprehension guarantees the existence of, say,
the set of all pairs $(p,q)$ such that $\uw p \supseteq \uw q$.
Therefore we showed that \ACA\ proves the equivalence of the six
characterizations of interval orders for countable partial orders.

Our goal is to obtain sharper results, in particular showing that
all equivalences can be proved in \WKL\ (which is strictly weaker
than \ACA). We will in fact use the results of this section about
finite partial orders to prove results about infinite partial orders
without resorting to the full power of \ACA.
\end{remark}

The following fact about the conjoint linear quasi-order will be
useful in the proof of Theorem \ref{tpt->io}.

\begin{lemma}\label{+-}
Let $\PP^*$ be the conjoint linear quasi-order associated to the
finite partial order \PP\ and let $p \in P$. Then:
\begin{itemize}
  \item
either $p^-$ is a minimum in $\PP^*$ (i.e.\ $\forall x \in P^*\, p^-
\leq^* x$) or there exists $q \in P$, $q \neq p$, such that $q^+$ is
an immediate predecessor of $p^-$ in $\PP^*$ (i.e.\ $x <^* p^-$
implies $x \leq^* q^+$ for all $x \in P^*$);
  \item
either $p^+$ is a maximum in $\PP^*$ (i.e.\ $\forall x \in P^*\, x
\leq^* p^+$) or there exists $q \in P$, $q \neq p$, such that $q^-$
is an immediate successor of $p^+$ in $\PP^*$ (i.e.\ $p^+ <^* x$
implies $q^- \leq^* x$ for all $x \in P^*$).
\end{itemize}
\end{lemma}
\begin{proof}
We prove the first statement (the second is proved similarly). Since
\PP\ and $\PP^*$ are finite, if $p^-$ is not minimal in $\PP^*$
there exists $x \in P^*$ which is an immediate predecessor of $p^-$.

To show that $x = q^+$ for some $q$, it suffices to show that for
every $r \in P$ with $r^- <^* p^-$ there exists $q \in P$ with $r^-
\leq^* q^+ <^* p^-$. Indeed, $r^- <^* p^-$ means $\dw r \subsetneqq
\dw p$ and there exists $q \in \dw p \setminus \dw r$. Then $q
\nless_P r$ and $q <_P p$ which imply $r^- \leq^* q^+$ and $q^+ <^*
p^-$.

It is obvious that $q \neq p$, since $p^- <^* p^+$.
\end{proof}

\section{Proofs in \RCA}\label{sect:RCA}
We start this section with the quite simple proof of the upper
upwards pointing arrow of Figure \ref{diagram} is provable in \RCA.

\begin{theorem}\label{clo->dist}
(\RCA) Every closed interval order is a distinguishing interval
order.
\end{theorem}
\begin{proof}
Let \PP\ be a closed interval order and let \LL, $f_0$ and $f_1$
witness this. Let $P^* = \set{p^+, p^-}{p \in P}$ and $L' = L \cup
P^*$ (we are assuming $L \cap P^* = \emptyset$).

We would like to define a linear order $\leq_{L'}$ on $L'$ so that
the maps $p \mapsto p^-$ and $p \mapsto p^+$ witness that \PP\ is a
distinguishing interval order. We first describe $\leq_{L'}$
informally: the restriction of $\leq_{L'}$ to $L$ coincides with
$\leq_L$, and $p^+$ and $p^-$ are placed respectively \lq\lq just
above $f_1(p)$\rq\rq\ and \lq\lq just below $f_0(p)$\rq\rq; if
distinct $p$ and $q$ are such that $f_1(p) = f_1(q)$ then $p^+$ and
$q^+$ are placed according to $\leq_\N$, and similarly for $p^-$ and
$q^-$ when $f_0(p) = f_0(q)$; if $f_0(p) = f_1(q)$ then $p^-$ is
below $q^+$.

To simplify the explicit definition of $\leq_{L'}$, we can exclude
the elements not belonging to the range of the functions we have in
mind, and therefore consider only the restriction of $\leq_{L'}$ to
$P^*$. Thus we set, for every $p, q \in P$:
\begin{align*}
p^+ \leq_{L'} q^+ & \iff f_1(p) <_L f_1(q) \lor (f_1(p) = f_1(q) \land p \leq_\N q);\\
p^- \leq_{L'} q^- & \iff f_0(p) <_L f_0(q) \lor (f_0(p) = f_0(q) \land p \leq_\N q);\\
p^+ \leq_{L'} q^- & \iff f_1(p) <_L f_0(q);\\
p^- \leq_{L'} q^+ & \iff f_0(p) \leq_L f_1(q).
\end{align*}
It is left to the reader checking that $\LL' = (P^*, {\leq_{L'}})$
is a linear order. We define $f'_0, f'_1: P \to P^*$ by $f'_0 (p) =
p^-$ and $f'_1 (p) = p^+$, and again we leave to the reader checking
that conditions (c1--2) and (c4) of definition \ref{def:io} hold.
Therefore \PP\ is a distinguishing interval order.
\end{proof}

We now show that also the bottom upwards pointing arrow of Figure
\ref{diagram} is provable in \RCA.

\begin{theorem}\label{tpt->io}
(\RCA) Every partial order not containing \tpt\ is an interval
order.
\end{theorem}
\begin{proof}
Let \PP\ be a partial order not containing \tpt. Let
$\set{p_n}{n>0}$ be an enumeration of $P$ (notice that for
notational convenience we start our enumeration from $1$). If $s \in
\N$ let $\PP_s = (\set{p_n}{0<n \leq s}, {\leq_P})$ and let
$\PP_s^*$ be the conjoint linear quasi-order associated to the
finite partial order $\PP_s$. We have $P^*_{s-1} \subset P^*_s$ and
we can investigate which relations are preserved from $\PP_{s-1}^*$
to $\PP_s^*$.

\begin{claim}\label{preservation}
$x <^*_{s-1} y$ implies $x <^*_s y$ for every $x,y \in P^*_{s-1}$.
\end{claim}
\begin{proof}
If exactly one of $x$ and $y$ is in $P_{s-1}^+$ (and the other is in
$P_{s-1}^-$) the claim follows immediately from the definition of
conjoint linear quasi-order. If $x,y \in P_{s-1}^+$, say $x = p_n^+$
and $y = p_m^+$, then $x <^*_{s-1} y$ means that $\uw
{p_n}^{\PP_{s-1}} \supsetneqq \uw {p_m}^{\PP_{s-1}}$.  Since $\uw
{p_i}^{\PP_s} \cap P_{s-1}^* = \uw {p_i}^{\PP_{s-1}}$, $\uw
{p_n}^{\PP_s} \subseteq \uw {p_m}^{\PP_s}$ cannot hold and, by Lemma
\ref{comparability} (which uses the hypothesis that \PP\ does not
contain \tpt), $\uw {p_n}^{\PP_s} \supsetneqq \uw {p_m}^{\PP_s}$,
i.e.\ $x <^*_s y$. The argument for the case $x,y \in P_{s-1}^-$ is
similar.
\end{proof}

On the other hand it is obvious that $x \equiv^*_{s-1} y$ does not
imply $x \equiv^*_s y$, e.g.\ if $x = p_n^+$, $y = p_m^+$, $\uw
{p_n}^{\PP_{s-1}} = \uw {p_m}^{\PP_{s-1}}$, $p_n <_\PP p_s$, and
$p_m \not <_\PP p_s$. We say that \emph{$x$ is separated below at
$s$} if for some $y$ we have $x \equiv^*_{s-1} y$ and $x <^*_s y$.
Analogously, \emph{$x$ is separated above at $s$} if for some $y$ we
have $x \equiv^*_{s-1} y$ and $y <^*_s x$.

\begin{claim}\label{separated}
At most one $\equiv^*_{s-1}$-equivalence class contained in
$P_{s-1}^+$ (recall Remark \ref{all1type}) contains elements
separated at $s$ (and the same for $\equiv^*_{s-1}$-equivalence
classes contained in $P_{s-1}^-$).
\end{claim}
\begin{proof}
Notice that by Lemma \ref{+-} $p_n^+$ can be separated at $s$ only
if $x <^*_s p_s^- <^*_s y$ for some $x,y \equiv^*_{s-1} p_n^+$. By
the previous claim, this can happen for the elements of at most one
$\equiv^*_{s-1}$-equivalence class.
\end{proof}

We define a linear quasi-order $\LL = (L, {\leq_L})$ where
\[
L = \set{x_n^k}{n \in \N \land n>0 \land k \in \Z \land n \leq |k|}.
\]
If $s \in \N$ let $L_s = \set{x_n^k \in L}{n \leq |k| \leq s}$. We
define $\leq_L$ by stages, so that at stage $s$ $\leq_L$ is defined
on the finite set $L_s$ and satisfies the following conditions:
\begin{enumerate}[(i)]
 \item
the set $\set{x_n^s, x_n^{-s}}{n \leq s} \subseteq L_s$ is ordered
by $\leq_L$ according to $\PP_s^*$, where $x_n^s$ and $x_n^{-s}$
replace respectively $p_n^+$ and $p_n^-$;
 \item
if $n<s$ then $x_n^{-s} <_L x_n^{-s+1}$ and $x_n^s >_L x_n^{s-1}$;
 \item
if $n<s$ and $y \in L_{s-1}$ then neither $x_n^{-s} \leq_L y <_L
x_n^{-s+1}$ nor $x_n^{s-1} <_L y  \leq_L x_n^s$ hold.
\end{enumerate}
An easy induction using (i) and (ii) yields $x_n^k <_L x_n^h$ if and
only if $k <_\Z h$. Notice also that (i) and (iii) imply $x_n^k
\not\equiv_L x_m^h$ whenever $k \neq h$.

Since $L_0 = \emptyset$ at stage $0$ there is nothing to do.

Let $s>0$ and suppose we have defined $\leq_L$ on $L_{s-1}$
satisfying (i--iii). To define $\leq_L$ on $L_s$ it suffices to
describe the position of the $x_n^s$'s and $x_n^{-s}$'s for $n \leq
s$.

First consider $x_n^s$ for $n<s$. If $p_n^+$ is not separated above
at $s$ then $x_n^s$ is an immediate successor (among the elements of
$L_s$) of $x_n^{s-1}$. If $p_n^+$ is separated above at $s$, fix
$p_m^+$ which is separated below at $s$. By Claim \ref{separated} we
have $p_m^+ \equiv^*_{s-1} p_n^+$, and hence by (i) $x_m^{s-1}
\equiv_L x_n^{s-1}$. Let $x_n^s$ be an immediate successor of
$x_s^{-s}$, which is an immediate successor of $x_m^s$ (which, by
the first clause of the present definition, is an immediate
successor of $x_m^{s-1} \equiv_L x_n^{s-1}$). The position of
$x_n^{-s}$ for $n<s$ is established similarly: if $p_n^-$ is not
separated below at $s$ then $x_n^{-s}$ is an immediate predecessor
of $x_n^{-s+1}$, otherwise fix $p_m^-$ which is separated above at
$s$ and let $x_n^{-s}$ be an immediate predecessor of $x_s^s$, which
is an immediate predecessor of $x_m^{-s}$.

If $p_s^+ \equiv_{P^*_s} p_n^+$ for some $n<s$ then set $x_s^s
\equiv_L x_n^s$, and similarly if $p_s^- \equiv_{P^*_s} p_n^-$ for
some $n<s$ set $x_s^{-s} \equiv_L x_n^{-s}$. If the previous case
does not hold and $p_s^+$ is the maximum in $\PP^*_s$ then $x_s^s$
is the maximum in $L_s$. Similarly if $p_s^-$ is the minimum in
$\PP^*_s$ then $x_s^{-s}$ is the minimum in $L_s$. If the position
of $x_s^s$ is not yet determined, by Lemma \ref{+-} $p_s^+$ is the
immediate predecessor in $\PP^*_s$ of some $p_n^-$ with $n<s$: let
$x_s^s$ be the immediate predecessor of $x_n^{-s}$ in $L_s$.
Similarly if $p_s^-$ is the immediate successor in $\PP^*_s$ of some
$p_n^+$ with $n<s$, let $x_s^{-s}$ be the immediate successor of
$x_n^s$ in $L_s$.

Notice that the latter part of the definition is compatible with the
positions of $x_s^s$ and $x_s^{-s}$ given earlier in some cases
(i.e.\ if some $p_n^-$ or $p_n^+$ is separated at $s$) above: in
fact if $p_m^+$ and $p_n^+$ are separated below and above,
respectively, at $s$ then $p_s^-$ is an immediate successor in
$\PP^*_s$ of $p_m^+$ (and similarly for the other case).

It is straightforward to check that $\leq_L$ restricted to $L_s$
satisfies (i--iii).

The definition of \LL\ is thus complete. We need to define $F
\subseteq P \times L$, and we would like to set
\[
F = \set{(p_n,x_m^k)}{\exists s\; x_n^{-s} \leq_L x_m^k \leq_L
x_n^s}.
\]
To show the existence of $F$ in \RCA, we need to prove that the
\SI01 formula appearing in the above definition is provably \DE01.

\begin{claim}\label{DE}
If $t = \max(|k|,n)$ then $\exists s\; x_n^{-s} \leq_L x_m^k \leq_L
x_n^s$ is equivalent to $x_n^{-t} \leq_L x_m^k \leq_L x_n^t$.
\end{claim}
\begin{proof}
One direction of the equivalence is obvious, so assume that
$x_n^{-s} \leq_L x_m^k \leq_L x_n^s$ for some $s \neq t$. If $s<t$
the conclusion follows immediately from $x_n^{-t} <_L x_n^{-s}$ and
$x_n^s <_L x_n^t$. If $s>t$ then $x_m^k \in L_{s-1}$ (because $m
\leq |k| \leq t < s$) and $n<s$: hence by (iii) we have $x_n^{-s+1}
\leq_L x_m^k \leq_L x_n^{s-1}$. Repeating this argument we obtain
$x_n^{-t} \leq_L x_m^k \leq_L x_n^t$.
\end{proof}

Claim \ref{DE} shows that $F$ exists. It is immediate that (i1) is
satisfied, so we need only to check (i2). If $p_n <_P p_m$ then by
(i) we have $x_n^s <_L x_m^{-s}$ for every $s \geq \max (n,m)$ and
this easily implies $\forall x \in F(p_n)\, \forall y \in F(p_m)\, x
<_L y$. If $p_n \nless_P p_m$ then $x_m^{-s} <_L x_n^s$ where $s =
\max (n,m)$: since $x_n^s \in F(p_n)$ and $x_m^{-s} \in F(p_m)$,
$\forall x \in F(p_n)\, \forall y \in F(p_m)\, x <_L y$
fails.\setcounter{claim}{0}
\end{proof}

\section{Equivalences with \WKL}\label{sect:WKL}
We first show that \WKL\ suffices to prove that the six
characterizations of interval orders we introduced are equivalent.

\begin{lemma}\label{forward}
(\WKL) Every partial order not containing \tpt\ is a distinguishing
interval order.
\end{lemma}
\begin{proof}
Let \PP\ be a partial order not containing \tpt. By Lemma
\ref{finite BCT} we can assume $P$ is infinite and let $\set{p_n}{n
\in \N}$ be a one-to-one enumeration of $P$. If $s \in \N$ let
$\PP_s = (\set{p_n}{n \leq s}, {\leq_P})$ and $\PP_s^*$ be the
conjoint linear quasi-order associated to the finite partial order
$\PP_s$. $\PP_s^*$ is a linear quasi-order by Lemma \ref{finite cio}
because \PP, and hence the finite partial order $\PP_s$, does not
contain \tpt.

Let $T$ be the set defined by setting $\sigma \in T$ if and only if
$\sigma$ is a finite sequence of length $\lh(\sigma)$ such that for
all $s,t < \lh(\sigma)$:
\begin{enumerate}
    \item
$\sigma(s)$ is (the code for) a linear order (denoted by
$\leq_{\sigma(s)}$) compatible with $\PP^*_s$ (see Definition
\ref{compatible});
    \item
if $s<t < \lh(\sigma)$ then $\sigma(t)$ extends $\sigma(s)$, i.e.\
$x \leq_{\sigma(s)} y \iff x \leq_{\sigma(t)} y$ for all $x,y \in
P^*_s$.
\end{enumerate}
$T$ exists by \DE01-comprehension. It is immediate that $T$ is a
tree. Since $\sigma(s)$ can assume only finitely many values
(corresponding to the (codes of the) finitely many linear orders on
the finite set $P^*_s$), $T$ is bounded in the sense of
\cite[Definition IV.1.3]{sosoa}. By Lemma \ref{existcompatible} for
every $s$ there exists a linear order compatible with $\PP^*_s$. By
taking its restrictions to $P^*_t$ for $t<s$ we construct a sequence
in $T$ of length $s$. Thus $T$ is infinite.

By Bounded K\"{o}nig's Lemma, which is provable in \WKL\ (\cite[Lemma
IV.1.4]{sosoa}), $T$ has an infinite path. This path is a sequence
$\set{\alpha(s)}{s \in \N}$ of (codes for) finite linear orders,
each one extending the previous ones and such that $\alpha(s)$ is
compatible with $\PP^*_s$. If $x,y \in P^*$ let $x \leq_L y$ if and
only if $x \leq_{\alpha(s)} y$ for any (or, equivalently, each) $s$
with $x,y \in P^*_s$. (Notice that here we are considering $P^*$
just as a set, without the ordering $\leq^*_\PP$ which is not
definable in \WKL.) $\leq_L$ exists by \DE01-comprehension.

It is straightforward to check that $(P^*,\leq_L)$ is a linear order
and that conditions (c1--2) and (c4) are satisfied by the functions
$p \mapsto p^-$, $p \mapsto p^+$ (because they are satisfied by each
$\leq_{\alpha(s)}$, by the proof of Lemma \ref{finite BCT}). Hence
\PP\ is a distinguishing interval order.
\end{proof}

\begin{cor}\label{WKL}
(\WKL) The five notions of interval order of Definition \ref{def:io}
and the property of not containing \tpt\ are all equivalent.
\end{cor}
\begin{proof}
This follows from Theorem \ref{obvious} and Lemma \ref{forward}.
\end{proof}

We now show that the implications that cannot be obtained by
composing arrows appearing in Figure \ref{diagram} are equivalent to
\WKL. In particular these implications are not provable in \RCA.

The following well-known characterization of \WKL\ (\cite[Lemma
IV.4.4]{sosoa}) is useful.

\begin{lemma}\label{functions}
(\RCA) The following are equivalent:
\begin{enumerate}[(i)]
  \item \WKL;
  \item
if $f,g: \N \to \N$ are one-to-one functions such that $\forall
n,m\; f(n) \neq g(m)$ then there exists a set $X$ such that $\forall
n (f(n) \in X \land g(n) \notin X)$.
\end{enumerate}
\end{lemma}

\begin{lemma}\label{rev1}
(\RCA) If every interval order is a 1-1 interval order then \WKL\
holds.
\end{lemma}
\begin{proof}
We will show that under our hypothesis (ii) of Lemma \ref{functions}
holds. Fix one-to-one functions $f,g: \N \to \N$ such that $\forall
n,m\; f(n) \neq g(m)$. We want to find a set $X$ such that $\forall
n (f(n) \in X \land g(n) \notin X)$.

We define a partial order $\leq_P$ on the set $P = \bigcup_{k \in
\N} P_k$, where $P_k = \{a_k,b_k\} \cup \set{c^n_k}{n \in \N}$ for
each $k$. If $p \in P_k$ and $q \in P_h$ with $k \neq h$ we set $p
\leq_P q$ if and only if $k<_\N h$. The elements of each $P_k$ are
pairwise $\leq_P$-incomparable with the following exceptions:
\begin{itemize}
 \item
if $n$ is such that $f(n)=k$ then $c^n_k <_P a_k <_P c^{n+1}_k$;
 \item
if $n$ is such that $g(n)=k$ then $c^n_k <_P b_k <_P c^{n+1}_k$.
\end{itemize}
Notice that our hypothesis on $f$ and $g$ imply that for each $k$ at
most one of the two possibilities occurs, and for at most one $n$.
$\leq_P$ can be defined within \RCA.

Let $\PP = (P, {\leq_P})$: it is immediate that \PP\ does not
contain \tpt. By Theorem \ref{tpt->io} \PP\ is an interval order and
by our hypothesis \PP\ is a 1-1 interval order. Hence there exist a
linear order $\LL = (L, {\leq_L})$ and $F \subseteq P \times L$
satisfying conditions (i1--3) of Definition \ref{def:io}. Let
$\varphi(k)$ and $\psi(k)$ be the \PI01 formulas
\[
F(a_k) \subseteq F(b_k) \qquad \text{and} \qquad F(b_k) \subseteq
F(a_k),
\]
respectively. Since (i3) holds (i.e.\ $F$ is one-to-one) we have
$\forall k\; \neg (\varphi(k) \land \psi(k))$ and we are in the
hypothesis of \PI01-separation (\cite[Exercise IV.4.8]{sosoa}),
which is provable in \RCA: hence there exists a set $X$ satisfying
\[
\forall k ((\varphi(k) \implies k \in X) \land (\psi(k) \implies k
\notin X)).
\]

We claim that $X$ satisfies also $\forall n (f(n) \in X \land g(n)
\notin X)$, thus completing the proof. To this end it suffices to
show that $\exists n\; f(n)=k$ implies $\varphi(k)$ and $\exists n\;
g(n)=k$ implies $\psi(k)$.

We prove only the first of these implications, the second being
similar. Suppose $n$ is such that $f(n)=k$: then $c^n_k <_P a_k <_P
c^{n+1}_k$, $c^n_k \nleq_P b_k$, and $b_k \nleq_P c^{n+1}_k$. The
last two conditions and (i2) imply the existence of $x \in
F(c^n_k)$, $x' \in F(b_k)$, $y \in F(c^{n+1}_k)$, and $y' \in
F(b_k)$ such that $x' \leq_L x$ and $y \leq_L y'$. By the first
condition and (i2), for all $z \in F(a_k)$ we have $x <_L z <_L y$,
and hence $x' <_L z <_L y'$. Now we use (i1), obtaining $F(a_k)
\subseteq F(b_k)$, i.e.\ $\varphi(k)$.
\end{proof}

\begin{lemma}\label{rev2}
(\RCA) If every 1-1 interval order is a closed interval order then
\WKL\ holds.
\end{lemma}
\begin{proof}
Again we will show that under our hypothesis (ii) of Lemma
\ref{functions} holds and we fix one-to-one functions $f,g: \N \to
\N$ such that $\forall n,m\; f(n) \neq g(m)$. We want to find $X$
such that $\forall n (f(n) \in X \land g(n) \notin X)$.

We define a partial order $\leq_P$ on the set $P = \bigcup_{k \in
\N} P_k$, where $P_k = \{a_k,b_k,c_k\} \cup \set{d^n_k}{n \in \N}$
for each $k$. As in the previous proof, if $p \in P_k$ and $q \in
P_h$ with $k \neq h$ we set $p \leq_P q$ if and only if $k<_\N h$.
Within each $P_k$ we have:
\begin{itemize}
  \item
$a_k \perp_P b_k$, $a_k \perp_P c_k$, and $c_k <_P b_k$;
  \item
if $f(n) \neq k \neq f(m)$ and $g(n) \neq k \neq g(m)$ then $d^n_k
<_P d^m_k$ if and only if $n<_\N m$;
  \item
if $f(n) \neq k$ and $g(n) \neq k$ then $a_k, b_k, c_k <_P d^n_k$;
  \item
if $f(n) = k$ and $m \neq n$ then $a_k, c_k <_P d^n_k <_P d^m_k$ and
$b_k \perp_P d^n_k$;
  \item
if $g(n) = k$ and $m \neq n$ then $b_k, c_k <_P d^n_k <_P d^m_k$ and
$a_k \perp_P d^n_k$.
\end{itemize}
\begin{figure}
\includegraphics{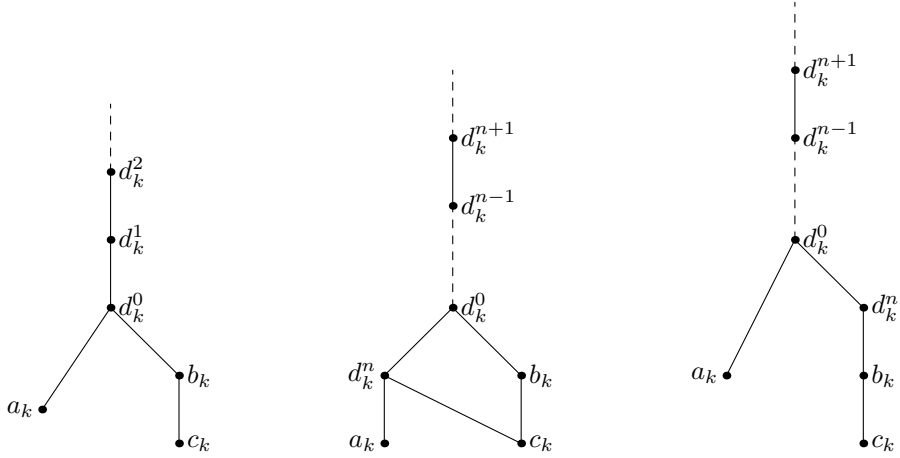}
\caption{\label{fig:cases}The three cases of $\leq_P$ restricted to
$P_k$ in the proof of Lemma \ref{rev2}: from left to right $\forall
n\; f(n) \neq k \neq g(n)$, $f(n)=k$, and $g(n)=k$.}
\end{figure}
Figure \ref{fig:cases} contains the Hasse diagram of the most
significant part of the restriction of $\leq_P$ to $P_k$ in the
three possible cases.

$\leq_P$ can be defined in \RCA. Let $\PP = (P, {\leq_P})$.

\begin{claim}
\PP\ is a 1-1 interval order.
\end{claim}
\begin{proof}
It is easy to check that \PP\ does not contain \tpt\ and hence it is
an interval order by Theorem \ref{tpt->io}. By Lemma
\ref{1-1condition} to prove the claim it suffices to show that
\[
\forall p,q \in P\; (p \neq q \implies \uw p \neq \uw q \lor \dw p
\neq \dw q).
\]

Fix $p,q \in P$ with $p \neq q$. If $p <_P q$ or $q <_P p$ then both
$\uw p \neq \uw q$ and $\dw p \neq \dw q$ hold. If $p \perp_P q$
then $p,q \in P_k$ for some $k$, and we consider the different
possibilities. In each case we exhibit an element of $P$ witnessing
either $\uw p \neq \uw q$ or $\dw p \neq \dw q$: $c_k \in \dw {b_k}
\setminus \dw {a_k}$, $b_k \in \uw {c_k} \setminus \uw {a_k}$, if
$f(n)=k$ then $a_k \in \dw {d_k^n} \setminus \dw {b_k}$, and if
$g(n)=k$ then $b_k \in \dw {d_k^n} \setminus \dw {a_k}$.
\end{proof}

By our hypothesis \PP\ is a closed interval order and there exist a
linear order $\LL$ and $f_0,f_1: P \to L$ satisfying (c1-2). Let $X
= \set{k}{f_1(a_k) \leq_L f_1(b_k)}$. To complete the proof we need
to check that $f(n) \in X$ and $g(n) \notin X$ for every $n$. If
$k=f(n)$ then $f_1(a_k) <_L f_0(d_k^n) \leq_L f_1(b_k)$ and $k \in
X$. If $k=g(n)$ then $f_1(b_k) <_L f_0(d_k^n) \leq_L f_1(a_k)$ and
$k \notin X$.\setcounter{claim}{0}
\end{proof}

We summarize our results in the following theorem (a few more
implications equivalent to \WKL\ can be stated using the information
contained in Figure \ref{diagram}, Corollary \ref{WKL}, and Lemmas
\ref{rev1} and \ref{rev2}).

\begin{theorem}\label{equivalence}
(\RCA) The following are equivalent:
\begin{enumerate}[(i)]
  \item \WKL;
  \item
every partial order not containing \tpt\ is a 1-1 interval order;
  \item
every interval order is a 1-1 interval order;
  \item
every 1-1 interval order is a distinguishing interval order;
  \item
every 1-1 interval order is a closed interval order.
\end{enumerate}
\end{theorem}
\begin{proof}
The forward direction, i.e.\ the fact that (i) implies each of
(ii)--(v), is a consequence of Corollary \ref{WKL}.

The implication (ii) $\implies$ (iii) follows from Theorem
\ref{obvious}(iv). Lemma \ref{rev1} shows that (iii) implies (i).
The implication (iv) $\implies$ (v) is immediate by Theorem
\ref{obvious}. Lemma \ref{rev2} shows that (v) implies (i).
\end{proof}

\section{Proper interval orders}\label{sect:proper}
In this section we deal with proper interval orders. Throughout most
of the section we point out the changes needed in the definitions
and proofs of \S\ref{sect:def}--\ref{sect:WKL}. However, Theorem
\ref{11->clo-proper} is new, because its statement without \lq\lq
proper'' is false by Lemma \ref{rev2}. The proof of Lemma \ref{rev3}
is also new, because the interval order used in the proof of Lemma
\ref{rev1} is not proper.

We start with the definitions and elementary facts corresponding to
Section \ref{sect:def}.

\begin{definition}\label{def:pio}
(\RCA) A partial order \PP\ is a \emph{proper interval order} if
there exist a linear order \LL\ and a set $F \subseteq P \times L$
such that (i1--2) of Definition \ref{def:io} hold and
moreover:\begin{enumerate}[\quad({i}1)]\setcounter{enumi}3
 \item
$F(p) \subseteq F(q)$ implies $F(p) = F(q)$ for all $p,q \in P$.
\end{enumerate}
\PP\ is a \emph{proper 1-1 interval order} if (i3) of Definition
\ref{def:io} holds as well.

\PP\ is a \emph{proper closed interval order} if there exist a
linear order \LL\ and functions $f_0,f_1: P \to L$ such that (c1--2)
of Definition \ref{def:io} hold and
moreover:\begin{enumerate}[\quad(c1)]\setcounter{enumi}4
 \item
$f_0(p) <_L f_0 (q)$ if and only if $f_1(p) <_L f_1 (q)$ for all
$p,q \in P$.
\end{enumerate}
\PP\ is a \emph{proper 1-1 closed interval order} if (c3) of
Definition \ref{def:io} holds as well. \PP\ is a \emph{proper
distinguishing interval order} if beside (c1--2) and (c5) we have
also (c4).
\end{definition}

\begin{definition}
(\RCA) A partial order \PP\ \emph{does not contain \tpo} if
\[
\forall p_0, p_1, p_2, q \in P (p_0 <_P p_1 <_P p_2 \implies p_0
\leq_P q \lor q \leq_P p_2).
\]
\end{definition}

\begin{lemma}\label{paired}
(\RCA) If \PP\ does not contain \tpo\ then for every $p,q \in P$ we
have either $\dw p \subseteq \dw q$ or $\uw p \subseteq \uw q$.
\end{lemma}
\begin{proof}
Towards a contradiction assume that $\dw p \nsubseteq \dw q$ and
$\uw p \nsubseteq \uw q$. If $p_0 \in \dw p \setminus \dw q$ and
$p_2 \in \uw p \setminus \uw q$, then $p_0,p,p_2,q$ witness that
\PP\ contains \tpo.
\end{proof}

\begin{theorem}\label{obvious-proper}
(\RCA)\begin{enumerate}[(i)]
 \item
Every proper (distinguishing) (1-1) (closed) interval order is a
(distinguishing) (1-1) (closed) interval order.
 \item
Every proper distinguishing interval order is a proper 1-1 closed
interval order.
 \item
Every proper 1-1 (closed) interval order is a proper (closed)
interval order.
 \item
Every proper (1-1) closed interval order is a proper (1-1) interval
order.
 \item
Every proper interval order contains neither \tpt\ nor \tpo.
\end{enumerate}\end{theorem}
\begin{proof}
Statement (i) is immediate from the definitions. The statements in
(ii--iv) are proved exactly as the corresponding statements in
Theorem \ref{obvious}.

To prove (v) let \PP\ be a proper interval order: by (i) above \PP\
is an interval order and by Theorem \ref{obvious}(iv) \PP\ does not
contain \tpt.

To show that \PP\ does not contain \tpo\ let $L$ and $F$ witness
that \PP\ is a proper interval order, and suppose towards a
contradiction that $p_0, p_1, p_2, q \in P$ are such that $p_0 <_P
p_1 <_P p_2$, $p_0 \nleq_P q$ and $q \nleq_P p_2$. The second
condition implies the existence of $x,y \in L$ such that $x \in
F(p_0)$, $y \in F(q)$, and $y \leq_L x$. Similarly by the third
condition there exist $y',x'$ such that $y' \in F(q)$, $x' \in
F(p_2)$, and $x' \leq_L y'$. For every $z \in F(p_1)$ the first
condition implies $x <_L z <_L x'$: this implies on one hand $y,y'
\notin F(p)$, and on the other hand $y <_L z <_L y'$ and hence $z
\in F(q)$ by (i1), for all $z \in F(p_1)$. Therefore $F(p_1)
\subsetneqq F(q)$, contradicting condition (i4).
\end{proof}

We now analyze finite partial orders containing neither \tpt\ not
\tpo, imitating what we did in Section \ref{sect:finite}.

\begin{definition}\label{def:properconjoint}
(\RCA) Given a finite partial order \PP\ let $P^\# = P^*$ be defined
as in Definition \ref{def:conjoint}. Define a binary relation
${\leq_\PP^\#}$ on $P^\#$ as follows:
\begin{align*}
p^+ \leq_\PP^\# q^+ & \iff \uw p^\PP \supsetneqq \uw q^\PP \lor
    (\uw p^\PP = \uw q^\PP \land \dw p^\PP \subseteq \dw q^\PP);\\
p^- \leq_\PP^\# q^- & \iff \dw p^\PP \subsetneqq \dw q^\PP \lor
    (\dw p^\PP = \dw q^\PP \land \uw p^\PP \supseteq \uw q^\PP);\\
p^+ \leq_\PP^\# q^- & \iff p <_P q;\\
p^- \leq_\PP^\# q^+ & \iff q \nless_P p.
\end{align*}
$\PP^\# = (P^\#, {\leq_\PP^\#})$ is the \emph{proper conjoint linear
quasi-order associated to \PP}. When \PP\ is clear from the context
we write $\leq^\#$ in place of $\leq_\PP^\#$.
\end{definition}

\begin{remark}
Notice that $\leq_\PP^\#$ and $\leq_\PP^*$ are defined on the same
set. It is immediate that ${\leq_\PP^\#} \subseteq {\leq_\PP^*}$,
and in general equality does not hold: in fact if $\uw p^\PP = \uw
q^\PP$ it is always the case that $p^+ \leq_\PP^* q^+$, while $p^+
\leq_\PP^\# q^+$ fails when $\dw p^\PP \nsubseteq \dw q^\PP$.
\end{remark}

The following lemma justifies the use of the words \lq\lq linear
quasi-order\rq\rq\ in Definition \ref{def:properconjoint}.

\begin{lemma}\label{finite pclq}
(\RCA) If \PP\ is a finite partial order which does not contain
\tpt\ then ${\leq^\#}$ is a linear quasi-order.

Moreover, if \PP\ does not contain \tpo\ then $\PP^\#$ and the
functions $p \mapsto p^-$, $p \mapsto p^+$ show that \PP\ is a
proper closed interval order.
\end{lemma}
\begin{proof}
The proofs that ${\leq^\#}$ is a linear quasi-order and that the
functions $p \mapsto p^-$, $p \mapsto p^+$ witness that \PP\ is a
closed interval order are identical to the same proofs for $\leq^*$
in Lemma \ref{finite cio}. Hence we need only to show that condition
(c5) of Definition \ref{def:pio} is met, i.e.\ that $p^- <^\# q^-$
if and only if $p^+ <^\# q^+$ for all $p,q \in P$.

Suppose $p,q \in P$ are such that $p^- <^\# q^-$ holds. Then either
$\dw p \subsetneqq \dw q$ or $\dw p = \dw q$ and $\uw p \supsetneqq
\uw q$. In the first case Lemma \ref{paired} implies that $\uw q
\subseteq \uw p$; even if $\uw q = \uw p$ we have $q^+ \nleq^\# p^+$
(because $\dw q \nsubseteq \dw p$) and hence $p^+ <^\# q^+$. In the
second case $p^+ <^\# q^+$ is immediate.

The reverse implication is proved similarly.
\end{proof}

\begin{remark}\label{all1type-proper}
Remark \ref{all1type} applies also to $\leq^\#$, i.e.\ each
$\equiv^\#$-equivalence class is contained in either $P^+$ or $P^-$.
Moreover $p^+ \equiv^\# q^+$ if and only if $\uw p^\PP = \uw q^\PP$
and $\dw p^\PP = \dw q^\PP$, if and only if $p^- \equiv^\# q^-$.
Therefore the $\equiv^\#$-equivalence classes contained in $P^+$ are
paired in a straightforward way with those contained in $P^-$.
\end{remark}

\begin{definition}\label{compatible-proper}
(\RCA) Given a finite partial order \PP\ which contains neither
\tpt\ nor \tpo, let $\PP^\#$ be the proper conjoint linear
quasi-order associated to \PP. A linear order $(P^\#, \leq_L)$ is
\emph{compatible with $\PP^\#$} if
\begin{gather*}
\forall x,y \in P^\# (x <^\# y \implies x <_L y),\\
\forall p,q \in P (p \neq q \land p^+ \equiv^\# q^+ \land p^+ <_L
q^+ \implies p^- <_L q^-), \quad \text{and}\\
\forall p,q \in P (p \neq q \land p^- \equiv^\# q^- \land p^- <_L
q^- \implies p^+ <_L q^+).
\end{gather*}
(Actually the second and third conditions imply each other.)
\end{definition}

\begin{remark}
Defining $(P^\#, \leq_L)$ compatible with $\PP^\#$ means defining a
linear order on each ${\equiv^\#}$-equivalence class, and keeping
the order between ${\equiv^\#}$-inequivalent elements unchanged.
Moreover we require that the linear orders on the
${\equiv^\#}$-equivalence classes containing $p^+$ and $p^-$ are the
same.
\end{remark}

\begin{lemma}\label{existcompatible-proper}
(\RCA) If \PP\ is a finite partial order which contains neither
\tpt\ nor \tpo\ then there exists a linear order compatible with
$\PP^\#$.
\end{lemma}
\begin{proof}
For example let
\[
x \leq_L y \iff x <^\# y \lor (x \equiv^\# y \land x \leq_\N y).
\]
$\leq_L$ is a linear order compatible with $\PP^\#$.
\end{proof}

\begin{lemma}\label{finite BCT-proper}
(\RCA) Any finite partial order which contains neither \tpt\ nor
\tpo\ is a proper distinguishing interval order.
\end{lemma}
\begin{proof}
Let \PP\ be a finite partial order which contains neither \tpt\ nor
\tpo, and, by Lemma \ref{existcompatible-proper}, $\leq_L$ a linear
order compatible with $\PP^\#$. Then $(P^\#, \leq_L)$ and the
functions $p \mapsto p^-$, $p \mapsto p^+$ show that \PP\ is a
proper distinguishing interval order. Indeed if $p \neq q$ and, say,
$p^+ \equiv^\# q^+$ then we have also $p^- \equiv^\# q^-$: if $p^+
<_L q^+$ then the second condition of Definition
\ref{compatible-proper} implies $p^- <_L q^-$.
\end{proof}

Combining Lemma \ref{finite BCT-proper} with Theorem
\ref{obvious-proper} we obtain that \RCA\ proves the equivalence of
the six characterizations of proper interval orders in the finite
case.

\begin{remark}\label{ACA-proper}
Remark \ref{ACA} applies also to what we have done with $\leq^\#$ in
the previous Lemmas, and we can conclude that \ACA\ suffices to
prove the equivalence of the six characterizations of proper
interval orders for countable partial orders.

As with interval orders, we will obtain sharper results also for
proper interval orders, in particular showing that all equivalences
can be proved in \WKL.
\end{remark}

\begin{remark}\label{not+-}
Notice that Lemma \ref{+-} does not hold with $\PP^\#$ in place of
$\PP^*$. If $P=\{p,q,r\}$ is ordered by $\leq_P$ as $\mathbf{2}
\oplus \mathbf{1}$ (i.e.\ the only nonreflexive relation is $p <_P
q$) then $p^- <^\# r^- <^\# p^+ <^\# q^- <^\# r^+ <^\# q^+$.
\end{remark}

Now we show that the upwards pointing implications of Figure
\ref{diagram-proper} are provable in \RCA, much as we did with
Figure \ref{diagram} in Section \ref{sect:RCA}.

\begin{theorem}\label{clo->dist-proper}
(\RCA) Every proper closed interval order is a proper distinguishing
interval order.
\end{theorem}
\begin{proof}
We can repeat the proof of Theorem \ref{clo->dist}. One needs only
to check that the construction preserves properness. We leave this
to the reader.
\end{proof}

As already noticed, the next Theorem has no counterpart for
arbitrary interval orders.

\begin{theorem}\label{11->clo-proper}
(\RCA) Every proper 1-1 interval order is a proper closed interval
order.
\end{theorem}
\begin{proof}
Let $\LL = (L, F)$ witness that the partial order \PP\ is a proper
1-1 interval order.

\begin{claim}
For all $p,q \in P$ the following are equivalent:
\begin{enumerate}[(1)]
 \item
$p = q \lor \exists x, y \in L (x \in F(p) \setminus F(q) \land y
\in F(q) \land x <_L y)$;
 \item
$\forall x, y \in L (x \in F(p) \setminus F(q) \land y \in F(q)
\implies x <_L y)$.
\end{enumerate}\end{claim}
\begin{proof}
First assume that (1) holds and (2) fails. Since $p=q$ implies (2),
there exist $x, y, x', y' \in L$ with $x,x' \in F(p) \setminus
F(q)$, $y,y' \in F(q)$, $x <_L y$ and $y' <_L x'$. Let $z \in F(q)$:
we have neither $z \leq_L x$ (because $x \notin F(q)$) nor $x'
\leq_L z$ (because $x' \notin F(q)$). Hence $x <_L z <_L x'$ and
$F(q) \subseteq F(p)$. Since it is immediate that $F(q) \neq F(p)$,
we are contradicting condition (i4) in definition \ref{def:pio}.

Now assume (2) holds and (1) fails, so that in particular $p \neq q$
and hence $F(p) \neq F(q)$ because condition (i3) holds. If $ F(p)
\setminus F(q) = \emptyset$ then $F(q) \subseteq F(p)$ and we are
again contradicting (i4). Therefore we can choose $x \in F(p)
\setminus F(q)$ and $y \in F(q)$: (2) implies $x <_L y$ and then we
have (1), against our assumption.
\end{proof}

Obviously (1) is \SI01 and (2) is \PI01. We denote either of them by
$\varphi(p,q)$: $\varphi$ is a provably \DE01 formula and we can use
it in the comprehension scheme. The following two claims about
$\varphi$ are useful.

\begin{claim}\label{claimuw}
$\varphi(p,q)$ implies $\uw q \subseteq \uw p$ and $\dw p \subseteq
\dw q$.
\end{claim}
\begin{proof}
Let $r \in \uw q$: to show $r \in \uw p$, i.e.\ $p <_P r$, by
(i2) it suffices to show that $x <_L z$ for all $x \in F(p)$ and $z
\in F(r)$. If $x \in F(q)$ this follows from $q <_P r$. If $x \in
F(p) \setminus F(q)$ let $y \in F(q)$: we have $x <_L y <_L z$ and
we are done.

The proof that $\dw p \subseteq \dw q$ is even simpler.
\end{proof}

\begin{claim}\label{claimcomp}
For every $p,q \in P$ either $\varphi(p,q)$ or $\varphi(q,p)$ holds.
\end{claim}
\begin{proof}
When $p=q$ the claim is obvious, so we assume $p \neq q$. Then $F(p)
\neq F(q)$ by (i3) and by (i4) $F(p) \setminus F(q)$ and $F(q)
\setminus F(p)$ are both nonempty. Let $x \in F(p) \setminus F(q)$
and $y \in F(q) \setminus F(p)$: if $x <_L y$ then $\varphi(p,q)$
holds, if $y <_L x$ then we have $\varphi(q,p)$.
\end{proof}

Let $P^\# = P^+ \cup P^-$ and define $\leq_{L'}$ by
\begin{align*}
p^+ \leq_{L'} q^+ & \iff \varphi(p,q);\\
p^- \leq_{L'} q^- & \iff \varphi(p,q);\\
p^+ \leq_{L'} q^- & \iff p <_P q;\\
p^- \leq_{L'} q^+ & \iff q \nless_P p.
\end{align*}
Reflexivity of $\leq_{L'}$ is immediate from the fact that
$\varphi(p,p)$ holds for every $p$. To check transitivity start by
noticing that using (2) it is immediate that $\varphi(p,q)$ and
$\varphi(q,r)$ imply $\varphi(p,r)$. This gives two of the eight
cases. The other four cases where some hypothesis is of the form
$\varphi(p,q)$, are easily handled using Claim \ref{claimuw}. Only
two cases are left:\begin{itemize}
 \item
if $p^+ \leq_{L'} q^- \leq_{L'} r^+$ then $p <_P q$ and $r \nless_P
q$. Thus there exist $z \in F(r)$ and $y \in F(q)$ with $y \leq_L
z$. Since $p \neq r$ we can pick $x \in F(p) \setminus F(r)$: we
have $x <_L y$ and hence $x <_L z$. Therefore $\varphi(p,r)$ and
$p^+ \leq_{L'} r^+$;
 \item
if $p^- \leq_{L'} q^+ \leq_{L'} r^-$ then $q \nless_P p$ and $q <_P
r$. Let $x \in F(p)$ and $y \in F(q)$ be such that $x \leq_L y$.
Since $p \neq r$ we can choose $z \in F(r) \setminus F(p)$: $x <_L
z$ follows immediately and hence we have that $\varphi(r,p)$ does
not hold. By Claim \ref{claimcomp} we have $\varphi(p,r)$ and $p^-
\leq_{L'} r^-$.
\end{itemize}
The fact that $(P^\#, {\leq_{L'}})$ is linear follows immediately
from the definition and Claim \ref{claimcomp}.

Define $f_0, f_1: P \to P^\#$ as usual by $f_0(p) = p^-$ and $f_1(p)
= p^+$. Conditions (c1--2) and (c5) follow immediately from the
definition of $\leq_{L'}$. Therefore \PP\ is a proper closed
interval order.\setcounter{claim}{0}
\end{proof}

\begin{remark}
The reader may have noticed the construction of the proof of Theorem
\ref{11->clo-proper} satisfies also condition (c4). Therefore the
proof actually shows that \RCA\ suffices to prove that every proper
1-1 interval order is a proper distinguishing interval order. This
result is also obtained combining the statements of Theorems
\ref{11->clo-proper} and \ref{clo->dist-proper}.
\end{remark}

\begin{theorem}\label{tpt+tpo->pio}
(\RCA) Every partial order which contains neither \tpt\ nor \tpo\ is
a proper interval order.
\end{theorem}
\begin{proof}
The proof follows the pattern of the proof of Theorem \ref{tpt->io}:
throughout the proof we replace $\PP^*_s$ with $\PP^\#_s$, the
proper conjoint linear quasi-order associated to $\PP_s$. We point
out only the spots where differences occur.

To prove the analogous of Claim \ref{preservation} we need to
consider the case of $n,m<s$ such that $p_n^+ <^\#_{s-1} p_m^+$
because $\uw {p_n}^{\PP_{s-1}} = \uw {p_m}^{\PP_{s-1}}$ and $\dw
{p_n}^{\PP_{s-1}} \subsetneqq \dw {p_m}^{\PP_{s-1}}$. Beside Lemma
\ref{comparability}, also Lemma \ref{paired} (which uses the
hypothesis that \PP\ does not contain \tpo) is needed here: since
$\dw {p_m}^{\PP_s} \nsubseteq \dw {p_n}^{\PP_s}$ we have $\uw
{p_m}^{\PP_s} \subseteq \uw {p_n}^{\PP_s}$ and therefore $\uw
{p_m}^{\PP_s} \supsetneqq \uw {p_n}^{\PP_s}$ cannot occur. Hence
$p_n^+ <^\#_s p_m^+$.

The analogous of Claim \ref{separated} states that at most two
$\equiv^\#_{s-1}$-equivalence class contained in $P_{s-1}^+$ contain
elements separated at $s$, and the same for
$\equiv^\#_{s-1}$-equivalence classes contained in $P_{s-1}^-$.

The definition of $\leq_L$ on $L_s$ requires considering a few more
possible situations. When $n<s$ and $p_n^+$ is separated above at
$s$, fix $p_m^+$ separated below at $s$ with $p_m^+ \equiv^\#_{s-1}
p_n^+$ and hence $x_m^{s-1} \equiv_L x_n^{s-1}$. If
$\uw{p_n}^{\PP_s} \subsetneqq \uw{p_m}^{\PP_s}$ then no changes are
needed, but now it might happen that $\uw{p_n}^{\PP_s} =
\uw{p_m}^{\PP_s}$ (because $\dw{p_n}^{\PP_s} \supsetneqq
\dw{p_m}^{\PP_s}$ forces $p_m^+ <^\#_s p_n^+$). In the latter case
$x_n^s$ is an immediate successor of $x_m^s$, which by the other
clauses in the definition is an immediate successor of $x_m^{s-1}
\equiv_L x_n^{s-1}$. If $p_n^-$ is separated below at $s$, act
similarly.

If $p_s^+$ is neither the maximum of $\PP^\#$ nor $\equiv^\#_s
p_n^+$ for some $n<s$ let $z \in P^\#_s$ be an immediate successor
of $p_n^+$ (now we cannot be sure that $z \in P_s^-$) and let
$x_s^s$ be an immediate predecessor of the element of $L_s \setminus
L_{s-1}$ which corresponds to $z$. Proceed analogously for
$x_s^{-s}$.

The definition of $F$ (including Claim \ref{DE}) and the proof that
\LL\ witnesses that \PP\ is an interval order needs no changes. Thus
we need only to show that condition (i4) is met. Assume $F(p_n)
\subseteq F(p_m)$ and fix $s \geq \max(n,m)$. By condition (iii) we
have $x_n^{-s} \leq_L x_m^{-s} <_L x_m^s \leq_L x_n^s$, and hence
$p_n^- \leq^\#_s p_m^- <^\#_s p_m^+ \leq^\#_s p_n^+$. By Lemma
\ref{finite pclq} this implies that $p_n^- \equiv^\#_s p_m^-$ and
$p_m^+ \equiv^\#_s p_n^+$, and hence $x_n^{-s} \equiv_L x_m^{-s}$
and $x_m^s \equiv_L x_n^s$. From the definition of $F$ we get
$F(p_n) = F(p_m)$, and the proof is complete.
\end{proof}

We now conclude with results similar to the one obtained in Section
\ref{sect:WKL}, showing that the implications missing from Figure
\ref{diagram-proper} are equivalent to \WKL.

\begin{lemma}\label{forward-proper}
(\WKL) Every partial order containing neither \tpt\ nor \tpo\ is a
proper distinguishing interval order.
\end{lemma}
\begin{proof}
The proof of Lemma \ref{forward} works without major changes,
replacing $\PP^*_s$ with $\PP^\#_s$. Obviously we use Lemmas
\ref{finite pclq}, \ref{existcompatible-proper}, and \ref{finite
BCT-proper} in place of Lemmas \ref{finite cio},
\ref{existcompatible}, and \ref{finite BCT}. Notice that since (c5)
is satisfied by each $\leq_{\alpha(s)}$ it is satisfied also by
$(P^\#, {\leq_L})$.
\end{proof}

\begin{cor}\label{WKL-proper}
(\WKL) The five notions of proper interval order of Definition
\ref{def:pio} and the property of containing neither \tpt\ nor \tpo\
are all equivalent.
\end{cor}
\begin{proof}
This follows from Theorem \ref{obvious-proper} and Lemma
\ref{forward-proper}.
\end{proof}

\begin{lemma}\label{rev3}
(\RCA) If every closed interval order which is also a proper
interval order is a proper closed interval order then \WKL\ holds.
\end{lemma}
\begin{proof}
We will show that under our hypothesis (ii) of Lemma \ref{functions}
holds. Fix one-to-one functions $f,g: \N \to \N$ such that $\forall
n,m\; f(n) \neq g(m)$. We want to find a set $X$ such that $\forall
n (f(n) \in X \land g(n) \notin X)$.

We define a partial order $\leq_P$ on the set $P = \bigcup_{k \in
\N} P_k$, where $P_k = \{a_k,b_k\} \cup \set{c^n_k}{n \in \N}$ for
each $k$. If $p \in P_k$ and $q \in P_h$ with $k \neq h$ we set $p
\leq_P q$ if and only if $k<_\N h$. The elements of each $P_k$ are
pairwise $\leq_P$-incomparable with the following exceptions:
\begin{itemize}
 \item
if $n$ is such that $f(n)=k$ then $a_k <_P c^n_k$;
 \item
if $n$ is such that $g(n)=k$ then $c^n_k <_P a_k$.
\end{itemize}
$\leq_P$ can be defined within \RCA. Let $\PP = (P, {\leq_P})$.

\begin{claim}\label{Pclosed}
\PP\ is a closed interval order.
\end{claim}
\begin{proof}
Let $\mathbf{N} = (\N, {\leq_\N})$ and define $f_0, f_1: \N \to P$
by setting
\begin{alignat*}2
f_0(a_k) = f_1(a_k) & = 3k+1;\\
f_0(b_k) & = 3k;\\
f_1(b_k) & = 3k +2;\\
f_0(c^n_k) & = 3k && \text{if $f(n) \neq k$};\\
f_1(c^n_k) & = 3k+2 \qquad && \text{if $g(n) \neq k$};\\
f_0(c^n_k) & = 3k+2 && \text{if $f(n) = k$};\\
f_1(c^n_k) & = 3k \qquad && \text{if $g(n) =k$}.
\end{alignat*}
It is straightforward to check that conditions (c1--2) of Definition
\ref{def:io} are met.
\end{proof}

\begin{claim}\label{Pproper}
\PP\ is a proper interval order.
\end{claim}
\begin{proof}
Claim \ref{Pclosed} and Theorem \ref{obvious} imply that \PP\ does
not contain \tpt. Our hypothesis on $f$ and $g$ imply that $c^n_k
<_P a_k <_P c^m_k$ cannot occur: hence \PP\ does not contain \tpo.
By Theorem \ref{tpt+tpo->pio}, \PP\ is a proper interval order.
\end{proof}

Claims \ref{Pclosed} and \ref{Pproper} and our hypothesis imply that
\PP\ is a proper closed interval order. Hence there exist a linear
order $\LL = (L, {\leq_L})$ and $f_0,f_1: P \to L$ satisfying
conditions (c1--2) of Definition \ref{def:io} and condition (c4) of
Definition \ref{def:pio}. Let $X = \set{k \in \N}{f_1(a_k) <_L
f_1(b_k)}$.

We now show that $X$ satisfies $\forall n (f(n) \in X \land g(n)
\notin X)$, thus completing the proof. If $f(n) = k$ then $a_k <_P
c^n_k$ and $b_k \nless_P c^n_k$: hence $f_1(a_k) <_L f_0(c^n_k)
\leq_L f_1(b_k)$ and $k \in X$. If $g(n) = k$ then $c^n_k <_P a_k$
and $c^n_k \nless_P b_k$: hence $f_0(b_k) \leq_L f_1(c^n_k) <_L
f_0(a_k)$. From $f_0(b_k) <_L f_0(a_k)$, (c4) yields $f_1(b_k) <_L
f_1(a_k)$ and hence $k \notin X$.
\end{proof}

\begin{theorem}\label{equivalence-proper}
(\RCA) The following are equivalent:
\begin{enumerate}[(i)]
  \item \WKL;
  \item
every partial order containing neither \tpt\ nor \tpo\ is a proper
1-1 interval order;
  \item
every partial order containing neither \tpt\ nor \tpo\ is a proper
closed interval order;
  \item
every proper interval order is a proper 1-1 interval order;
  \item
every closed interval order which is also a proper interval order is
a proper closed interval order.
\end{enumerate}
\end{theorem}
\begin{proof}
The forward direction, i.e.\ the fact that (i) implies each of
(ii)--(v), is a consequence of Corollary \ref{WKL-proper}.

The implications (ii) $\implies$ (iii) and (iv) $\implies$ (v)
follow from Theorem \ref{11->clo-proper}. Theorem
\ref{obvious-proper}(v) shows (ii) $\implies$ (iv). The implication
(iii) $\implies$ (v) is immediate by Theorem \ref{obvious-proper}.
Lemma \ref{rev3} shows that (v) implies (i).
\end{proof}

\bibliography{interval}
\bibliographystyle{hplain}

\end{document}